\documentclass[]{amsart} 
\usepackage[]{amsmath, amsthm, amsfonts}
\usepackage[all]{xy} 

\newcommand {\C} {{\mathbb C}}
 \newcommand {\R} {{\mathbb R}} 
\newcommand {\Z} {{\mathbb Z}}
 \newcommand {\Q} {{\mathbb Q}} 
\newcommand {\PP} {{\mathbb P}}
  
\newcommand {\F} {{\mathcal F}}

\newcommand {\E} {{\mathcal E}} 
\newcommand {\dt} {{\bullet}}

\newcommand {\LL} {{\mathcal L}}

\newcommand {\V} {{\mathcal V}}

 \newtheorem{thm}[subsection]{Theorem}
 \newtheorem{cor}[subsection]{Corollary}
 \newtheorem{lemma}[subsection]{Lemma}
 \newtheorem{prop}[subsection]{Proposition}

\begin{document}

 \title{Motivation for Hodge cycles}

 \author{ Donu Arapura} \thanks{Author partially supported by the NSF}
 \address{Department of Mathematics\\
   Purdue University\\
   West Lafayette, IN 47907\\
   U.S.A.}
 
 \maketitle

 Given two smooth projective varieties $X$ and $Y$ over a field, we
 say that $X$ {\em motivates } $Y$ or that $Y$ is {\em motivated by}
 $X$ if the motive of $Y$ is contained in the category generated from
 $X$ by taking sums, summands and products.  This notion has appeared
 implicitly in many places, but it seems useful to isolate it so as to
 state the following principle (lemma \ref{lemma:motivatinglemma}): if
 the Hodge (generalized Hodge, Lefschetz standard...) conjecture holds
 for $X$ and all its powers, then it holds for any variety motivated
 by it.  For the precise statement we can use homological motives,
 however, we find it more convenient to use the construction of
 motives due to Andr\'e \cite{andre} which has the advantage of
 yielding a (provably) semisimple Abelian category through which
 cohomology factors.

 Given a smooth complex projective variety $X$, we can take the
 dimension of the smallest variety that motivates $X$ as a measure of
 its complexity.  This number can be seen to be maximal for general
 varieties using work of Schoen\cite{schoen}; however, there are a
 number of interesting examples, discussed below, where this is small.
 Varieties motivated by curves are the simplest. For such varieties, a
 weak form of the Hodge conjecture, that Hodge cycles are motivated in
 Andr\'e's sense, holds unconditionally.  Such cycles are absolutely
 Hodge in Deligne's sense.  Next in line are varieties motivated by
 curves or surfaces. For these varieties we check that the Lefschetz
 standard conjecture of Grothendieck holds.

 There are a number of natural examples of varieties motivated by
 curves and surfaces.  These include Abelian varieties, uniruled
 threefolds and unirational fourfolds. These are checked by direct
 geometric arguments. For Abelian varieties, we first observe that the
 Jacobian of a curve is motivated by the curve. This generalizes to
 other moduli spaces.  We show that the moduli space of stable
 parabolic bundles over a curve is motivated by the curve, the Hilbert
 scheme of points over a surface is motivated by the surface, and
 likewise for the moduli space of stable vector bundles over an
 Abelian or K3 surface. For vector bundles over a curve, this result
 was first proved by del Ba\~no \cite{delb}.  For the Hilbert scheme
 of surface, this goes back to Cataldo-Migliorini \cite{cm}.  However,
 we give uniform and self contained treatments of these cases.  Using
 these results, we check the (generalized) Hodge conjecture for the
 above spaces in some cases, and the Lefschetz standard conjecture in
 all cases.

 This paper is essentially a refined and streamlined version of my
 preprint \cite{arapura}. My thanks to Y. Andr\'e and the referee
for making a number  of  useful suggestions.

 \section{Motives}

 Let $k$ be a field. Let $SPVar_k$ be the category of smooth
 projective (possibly reducible) varieties over $k$. The case of
 primary interest for us is $k=\C$.  Given an object of $SPVar_\C$,
 let $H^\dt(X)$ denote singular rational cohomology of $X^{an}$ with
 its canonical Hodge structure. This takes values in the category
 $PHS$ of finite direct sums of polarizable rational Hodge structures.
 The category $PHS$ is a semisimple $\Q$-linear Abelian category
 \cite[4.2.3]{deligne-hodge} with tensor products and duals.

 We call a full subcategory $\V$ of $SPVar_k$ {\em admissible} if it
 contains $spec\, k,\PP_k^1$ and is stable under products, disjoint
 unions, and connected components. Let $\langle X\rangle$ be the
 smallest admissible category containing a variety $X$. 

Given an admissible category $\V$ and an object $X\in
SPVar_k$,  Andr\'e \cite{andre}  has constructed
a graded  $\Q$-algbera $A^\dt_{mot}(X)$ called the
algebra of {\em motivated cycles} on $X$ modeled on $\V$.  
We refer likewise to elements of $A^\dt_{mot}(X_1\times X_2)$ as
motivated correspondences modeled on $\V$.
Fix a Weil cohomology $H^*(X)$, then we can regard $A^\dt_{mot}(X)$ as
a subalgebra of $H^{2*}(X)$. 
A class $\gamma\in A^\dt_{mot}(X)$ if and only
if there exists an object $Y\in \V$ and algebraic cycles 
$\alpha,\beta$ on $X\times Y$ such that 
$$\gamma = p_*(\alpha\cup *\beta),$$
where $p:X\times Y\to X$ is the projection, and $*$ is the 
Lefschetz involution with respect to
a product polarization \cite{andre}. Note that  $A^\dt_{mot}(X)$ contains
 the algebra of algebraic cycles on $X$, and it would coincide with it
assuming Grothendieck's standard conjectures. Motivated cycles  forms
a good replacement for algebraic cycles  in lieu of these conjectures.

 By an {\em
   intersection theory} on an admissible category $\V$, we mean a
 functor $R$ from $\V^{op}$ to commutative rings equipped with
 pushforwards satisfying the conditions of \cite[section 1]{manin}.
 There are several examples of interest to us:
 \begin{enumerate}
 \item $R=K_0$, the Grothendieck group of coherent sheaves.
 \item The quotient of the rationalized Chow ring $CH(\>)\otimes \Q$
   by an adequate equivalence relation (e.g. identity, homological, or
   numerical equivalence).
 \item The ring $R(\>)= A_{mot}(\>)$ of motivated cycles modeled on
   $\V$ as explained above. See \cite{andre}.
\end{enumerate}

In all but the first case $R$ has a grading.

Given the above data, we can form the category $Cor_R^u(\V)$ of
(ungraded) $R$-correspondences in $\V$ with the same objects as $\V$,
and $Hom(X,Y) = R(X\times Y)$.  Composition is given by
$$\beta\circ\alpha =p_{XZ*}(p_{XY}^*\alpha\cdot p_{YZ}^*\beta)$$
where
$p_{XZ}:X\times Y\times Z\to X\times Z,\ldots$ are the projections.
We write $Cor_R^u$ (respectively $Cor_R^u(X)$) etcetera for
$Cor_R^u(SPVar_k)$ (respectively $Cor_R^u(\langle X\rangle)$)
etcetera. In the cases, where $R$ has a grading, we define the
subcategory of graded correspondences $Cor_R\subset Cor_R^u$ by
restricting $$Hom_{Cor_R}(X,Y) = \prod_i R^{\dim X_i }(X_i\times Y)$$
where $X_i$ are the connected components of $X$. The category of
ungraded (respectively graded) $R$-motives $M_R^u(\V)$ ($M_R(\V)$) in
$\V$ is obtained by taking the pseudo-abelian completion of
$Cor_R^u(\V)$ ($Cor_R(\V)$) and inverting the so called Lefschetz
motive. Alternatively following \cite{jannsen2, scholl}, the objects
of $M_R(\V)$ can be regarded as triples $(X,p,m)$, with $X\in Ob\V$,
$p\in End(X)= R^{\dim X}(X\times X)$ an idempotent, and $m$ an integer
(we will also write this as $(X,p,0)(m)$). The morphisms are given by
$$Hom((X,p,m), (Y,q,n)) = q\circ [R^{\dim X-m+n}(X\times Y)]\circ p$$
When $R= {CH}(\>)\otimes \Q$ (respectively $R= A^*(\>) = im\,
CH^*(\>)\otimes \Q\to H^{2*}(\>)$), $ M_{CH}= M_R$ (respectively $
M_{hom}=M_R$) is called the category of Chow (homological) motives.
When $R=A_{mot}$ is the ring of motivated cycles, we call $M_A=M_R$
(respectively $M_A(\V)$) the category of Andr\'e motives (modeled on $\V$).

We have obvious functors $M_{CH}\to M_{hom}\to M_A$.  These categories
are all $\Q$-linear pseudo-Abelian categories with tensor products and
duals (see \cite{scholl}), and furthermore $M_A$ is semisimple Abelian
\cite{andre}.  We can associate a motive $[X]= (X,id,0)$ (in any of
the previous senses) to a variety $X\in \V$, and this yields a
contravariant functor by assigning to $f:X\to Y$ the transpose of its
graph.

Suppose $k=\C$.  Then the functors $H^\dt$ extend to {\em covariant}
functors on $M_A$ as follows.  First, recall that a correspondence
$\gamma\in Hom_{Cor_{A}}(X,Y)$ acts on cohomology by $\gamma_*(\alpha)
= p_{Y*}(p_X^*\alpha\cup [\gamma])$.  Given $(X,p,m)$ define
$$H^i(X,p,m) = p_*H^{i+2m}(X)(m) $$
where $(m)$ represents Tate twist
of the canonical Hodge structure.  If $$f\in Hom((X,p,m), (Y,q,n))$$
is given by $q\circ \gamma\circ p$, then $\gamma_*$ induces a morphism
of Hodge structures $$f_*:p_*H^{i+2m}(X)(m)\to q_*H^{i+2n}(Y)(n)$$
These rules yield a functor $H^i$ from $M_A$ into the category pure
polarizable Hodge structures weight $i$.  The functor $X\mapsto H(X) =
\oplus H^i(X)$ gives faithful additive embeddings of $M_{hom}$ and
$M_A$ into the $PHS$ (the faithfulness can be checked using Manin's
identity principle \cite{manin, scholl}).  Since $M_A$ is semisimple
Abelian, the additivity forces $H$ and $H^i$ to be exact on it as
well. $H$ also preserves tensor products and duals.  These Hodge
structures are not compatible with ungraded correspondences. However,
after adjusting weights and summing, the Hodge structures $$\tilde
H^{even}(X,p,m) = \bigoplus_j p_*H^{2j+m}(X)(j+m)$$
$$\tilde
H^{odd}(X,p,m) = \bigoplus_j p_*H^{2j+m+1}(X)(j+m)$$
will give
functors from $M_{A}^u\to PHS$. Furthermore, $X\mapsto \tilde
H(X)=\tilde H^{even}(X)\oplus\tilde H^{odd}(X)$ gives a faithful
embedding. When $k$ is arbitrary, similar remarks apply with $H$
replaced by $\ell$-adic cohomology.

For any admissible class $\V$, we can identify $M_A(\V)$ with a
subcategory of $M_A$. This need not be a full embedding, since the
notion of motivated cycles modeled on $\V$ may be more restrictive
than motivated cycles modeled on all of $SPVar_k$.  Let 
$M_A(\V)^{full}\subseteq M_A$
be the full subcategory generated by $M_A(\V)$. We say that a
smooth projective variety $Y$ is {\em motivated by} $\V$ (or a smooth
projective variety $X$) if $[Y]$ lies in $M_A(\V)$ (or $M_A(X)$). 
More precisely, this means that $[Y]$ is isomorphic in $M_A$
to an object of $M_A(\V)$.
Replacing $M_A(\V)$ by $M_A(\V)^{full}$ leads to the more flexible
(although harder to control) notion of weak motivation. 
For example, Andr\'e \cite{andre} has shown  that any K3 surface is 
weakly motivated by an Abelian variety. The corresponding result for
motivation  is unknown except in special cases such as for Kummer surfaces.
$X$ and $Y$ will be called (weakly) {\em co-motivated} if they are (weakly)
motivated by each other.

\begin{lemma}\label{lemma:subord}
  Given $X,Y$ in $SPVar_k$, $Y$ is (weakly) motivated by $X$ if and only if
  there exists a morphism $$f:\bigoplus_{m,n} [X]^{\otimes n}(m)\to
  [Y]$$
  in $M_A(X)$ ($M_A(X)^{full}$) inducing a surjection on cohomology.
\end{lemma}

\begin{proof}
  If $[Y]$ lies in $M_A(X)$, then it is a direct summand of some
  $\oplus [X]^{n}(m)$. Therefore projection yields the desired
  morphism $f$.

  Conversely, given a morphism $f$ as above. Since $H:M_A\to PHS$ is
  faithful and exact, it follows that $f$ is an epimorphism.
  Therefore $[Y]$ is a summand of $\oplus [X]^{n}(m)$, since $M_A$ is
  semisimple.
\end{proof}

\begin{cor}
  $Y$ is motivated by $X$ if there exists a surjective morphism of
  varieties $f:X^n\to Y$.
\end{cor}

\begin{proof}
  By taking general hyperplane sections, we can find a smooth
  $g:Z\hookrightarrow X^n$ such that $h=f\circ g$ is surjective, and
  $\dim Z = \dim Y$.  The map $h_*:H(Z)\to H(Y)$ is surjective since
  $\frac{1}{\deg h}h^*$ splits it. There $f_*$ is also surjective.
\end{proof}

\begin{lemma}\label{lemma:XvsAlbX}
  If $X$ is smooth projective variety, its Albanese $Alb(X)$ is
  motivated by $X$.  If $X$ is a smooth projective curve, $X$ and its
  Jacobian $J(X)$ are co-motivated.
\end{lemma}

\begin{proof}
  Let $\alpha:X\to Alb(X)$ be the Abel-Jacobi map. Since $Alb(X)$ is
  generated as a semigroup by the image of $\alpha$, the map $X^n\to
  Alb(X)$ given by $(x_1,\ldots x_n)\mapsto \alpha(x_1)+\ldots
  \alpha(x_n)$ is surjective for some $n$. This proves the first
  statement.

  Suppose that $X$ is a curve. We have just seen that $J(X)=Alb(X)$ is
  motivated by $X$. Since $\alpha^*$ induces a surjection on
  cohomology, $X$ is also motivated by $J(X)$.
\end{proof}

\begin{lemma}\label{lemma:ABtype}
  Suppose that $X$ and $Y$ are smooth projective varieties such that
  there exists a finite collection of motivated correspondences on
  $X\times Y$ modeled on $\langle X\rangle$ (respectively
$SPVar_k$) whose K\"unneth components, along $Y$, generate the
  cohomology ring $H(Y)$. Then $Y$ is motivated (respectively weakly
motivated) by $X$.
\end{lemma}

 \begin{proof}
   Let $d=\dim X$, and let $c_{i,j}\in A^{d+i}_{mot}(X\times Y)$
   denote the classes of the given correspondences. These induce
   morphisms $[X](-i)\to [Y]$ in $M_A$. Products $c_{i_1,j_1}\otimes
   \ldots c_{i_n,j_n}$ induce morphisms $$[X]^{\otimes
     n}(-i_1-i_2\ldots)\to [Y]^n\stackrel{\Delta^*}{\to} Y$$
   By
   assumption, a finite sum of these morphisms yield a map
   $$f:\bigoplus_{(i_1,\ldots i_n)} [X]^n(*)\to [Y]$$
   which induces a
   surjection on cohomology. Therefore we are done by lemma
   \ref{lemma:subord}.
 \end{proof}

\begin{lemma}\label{lemma:diagonal}
  Let $X$ and $Y$ be smooth projective varieties.  Suppose that the
  diagonal $\Delta\in H(Y\times Y)$ is contained in the algebra
  generated by products $\mu\times \mu'$, where $ \mu,\mu'\in H(Y)$
  are K\"unneth components of motivated correspondences on $X\times
  Y$   modeled on $\langle X\rangle$ (respectively
$SPVar_k$). Then $X$ motivates (respectively weakly motivates) $Y$.
\end{lemma}

\begin{proof}
  By assumption, we can write $\Delta= \sum n_{\nu\xi}\nu\times \xi$,
  where $\nu,\xi$ are products of K\"unneth components of
  correspondences on $X\times Y$.  Given $\alpha \in H(Y)$, we have
  $$\alpha= \Delta_*(\alpha) = \sum n_{\nu\xi}\langle
  \nu,\alpha\rangle \xi,$$
  where $\langle \nu,\alpha\rangle$ denotes
  $\int \nu\cup \alpha$.  Thus the hypothesis of lemma
  \ref{lemma:ABtype} is fulfilled.
\end{proof}

\begin{lemma}\label{lemma:HXgend}
  If $H^*(X)$ is generated as an algebra by elements of degree at most
  $d$, then $X$ is weakly motivated by a variety of dimension less than or
  equal to $d$.
\end{lemma}

\begin{proof}
  If $\dim X \le d$ there is nothing to prove. Otherwise, let $\iota':
  Z\hookrightarrow X$ be an intersection of $X$ with $\dim X -d$
  hyperplanes in general position.  We get an induced morphism
  $\iota:[X]\to [Z]$ of motives. Since $M_A$ is semisimple, there
  exists a morphism $\sigma:[Z]\to [X]$ satisfying $\iota\sigma \iota
  = \iota$ (this can be obtained as a composition of splittings
  $[Z]\to im(\iota)\to [X]$).  By the weak Lefschetz theorem,
  $\iota_*:H^{a}(X)\to H^a(Z)$ is injective for $a\le d$.  Therefore
  the map $H^a(Z)\to H^a(X)$ induced by $\sigma$ is surjective when
  $a\le d$.  By assumption, $H(X)$ is generated, as an algebra, by the
  elements in the images of these maps. Therefore we are done by lemma
  \ref{lemma:ABtype}.
\end{proof}

We can view the smallest $m$, for which $X$ is weakly motivated by an
$m$-fold, as a measure of the complexity of $X$.  For such an $m$, the
Hodge structures $H^i(X)$ would have to lie in the tensor category
generated by Hodge structures of level at most $m$.  Following Schoen
\cite{schoen}, we can find a Hodge theoretic obstruction to this.
Given a Hodge structure $H$, let $\mu(H)$ denote the level of the
induced Hodge structure on the Mumford-Tate Lie algebra of $H$.
(Recall that the level of a Hodge structure $G$ is $\max|p-q|$
such that $G^{pq}\not=0$.)
 We have that $\mu(H)$ is bounded above by the twice the level of $H$, and that
$$\mu(H_1\otimes H_2\otimes\ldots H_n) \le \max\{\mu(H_i)\}$$
\cite[pp
546-548]{schoen}. From this it follows that if $H$ lies in the tensor
category generated by Hodge structures of level at most $m$, then
$\mu(H)\le 2m$.  Let $\tau(X)$ be Schoen's invariant, which is half
the maximum of $\mu(H')$ as $H'$ varies over all irreducible Hodge
substructures $H^i(X)$ of level $i$ for all $i$.  Then from this
discussion, we find:

\begin{lemma}\label{lemma:schoen}
  If $X$ is weakly motivated by an $m$ dimensional variety, then $\tau(X)\le
  m$.
\end{lemma}

Schoen \cite{schoen} gives examples, such as general hypersurfaces of
large degree, where $\tau(X)=\dim X$

\section{Singular or non-projective varieties}

It will be convenient to extend the previous ideas to the category
$Var_\C$ of all varieties over $\C$.  If $X$ is a proper variety which
is a rational homology manifold, then it is as good as smooth for our
purposes.  In particular, we can attach a homological motive $[X]$ to it as
follows.  Since the rational cohomology of $X$ satisfies Poincar\'e
duality, we have a Gysin map $p_*$ for any resolution $p:\tilde X\to
X$.  We take $[X]= (\tilde X, p^*p_*,0)$, which is easily seen to be
well defined.

More general varieties give rise to mixed motives, in principle.
However we will only need the pure part of this structure. This is
analogous to passing from a mixed Hodge structure $H$ to the pure
structure $Gr_\dt^WH$.  Let $K^b(A)$ denote the homotopy category of
bounded complexes in an additive category $A$. This has a natural
triangulated structure.  When $A$ is Abelian, we have functors
$h^i:K^b(A)\to A$ given by taking $i$th cohomology. The following was
obtained by Gillet-Soul\'e{\cite{gillet}} and Guillen-Navarro
\cite{navarro}.

\begin{thm}\label{thm:gillet}
  Let $k=\C$, then for each $X\in ObVar$, there exists a well defined
  complex $W(X)\in ObK^b(M_{CH})$ such that
  \begin{enumerate}
  \item When $X$ is a smooth projective variety, $W(X) \cong [X]$.
  \item $W$ behaves contravariantly for proper maps.
  \item $W$ behaves covariantly for open immersions.
  \item $W(X\times Y)\cong W(X)\otimes W(Y)$.
  \item If $U\subset X$ is open, there is a natural distinguished
    triangle $$W(U)\to W(X)\to W(X-U)\to W(U)[1].$$
  \item $h^j(H^i(W(X))) = Gr_i^WH_c^{i+j}(X)$, where $H_c$ denotes
    cohomology with compact support.
  \end{enumerate}
\end{thm}

There are a few cases in which this complex can be made rather
explicit.  Given a divisor with normal crossings $D=\cup D_i$ on a
smooth projective variety $X$, then $W(X-D)$ can be realized by the
complex $$
[ X] \to \bigoplus_i [D_i]\to \bigoplus_{i,j} [D_i\cap
D_j]\ldots$$
with simplicial coboundaries \cite{gillet}. Item (6) is essentially
given in \cite[p 147]{gillet}, however it can be seen directly
for $X-D$ from the above complex.

 If $X$ is a
smooth projective variety with a finite group action such that the
quotient $X/G$ is a variety, $W(X/G)$ is isomorphic to $e[X]$ in
degree $0$, where $e=(1/\#G)\sum g\in \Q[G]$ \cite{banoN}.  In fact,
it is easy to combine these two cases to see that if $G$ acts on $(
X,D)$, $W(( X-D)/G)$ is given by $$
e[ X] \to \bigoplus_i e[D_i]\to
\ldots$$

Let $W_A(X)$ be the image of $W(X)$ in $K^b(M_A)$.  We write $Gr_j[X]$
for $h^{j}(W_A(X))$.  Under the embedding $H:M_A\to PHS$, $H(Gr_j[X])
= \oplus_i Gr_i^WH_c^{i+j}(X)$.  The discussion in the previous
paragraph implies that if $X$ is smooth with a smooth compactification
$\bar X$, $Gr_0[X]$ is a subobject of $[\bar X]$.

\begin{cor}
  If $U\subset Y$ is open, we have an exact sequence $$\ldots
  Gr_{j}[U]\to Gr_j[Y]\to Gr_j[Y-U]\ldots$$
  in $M_A$.
\end{cor}

We will say that an arbitrary variety $Y$ is (weakly) motivated by $\V$ if
$W_A(Y)$ is isomorphic a complex in $K^b(M_A(\V))$ (respectively
$K^b(M_A(\V)^{full})$). If $Y$ is smooth
and projective, these notions are equivalent to the previous definitions since
$W_A(Y)\cong [Y]$.

\begin{cor}
  If If $U\subset Y$ is open and any two of $U, Y,Y-U$ are (weakly) motivated
  by $\V$, then so is the third.  If $X$ and $Z$ are (weakly) motivated by
  $\V$, then so is $X\times Z$.
\end{cor}

\begin{proof}
  Since any vertex of a distinguished triangle can be constructed from
  the other two in terms of mapping cones, the first statement
  follows.  The second statement is evident from the theorem.
\end{proof}

By induction, we get:

\begin{cor}\label{cor:strat}
  If $Y$ is a smooth projective variety which can be expressed as a
  disjoint union $Y=\cup Y_i$ of locally closed varieties, such that
  $Y_i$ is (weakly) motivated by $\V$. Then so is $Y$.
\end{cor}

Let us say that a morphism $Y\to S$ is {\em cellular} if it is flat
and admits a decomposition $Y=\cup Y_i$ with $Y_i$ isomorphic to an
affine space fibration ${\mathbb A}^{n_i}_S$.

\begin{lemma}
  If $Y\to S$ is cellular, then $Y$ is motivated by $S$.
\end{lemma}

\begin{proof}
  This follows from the previous two corollaries.  Alternatively, it
  can be deduced from the isomorphism of graded Chow motives
  $$[Y]\cong \bigoplus [S](i),$$
  given in \cite[2.6]{scholl}.
\end{proof}

Combining this with the previous results gives.

\begin{cor}\label{cor:strat2}
  Suppose that $Y$ is a disjoint union $Y=\cup Y_i$ of subvarieties
  admitting cellular maps $Y_i\to Z_i$ with $Z_i$ (weakly) motivated by $\V$.
  Then $Y$ is (weakly) motivated by $\V$.
\end{cor}

\begin{cor}\label{cor:blowup}
  The blow up of a smooth projective variety $Y$ along a smooth center
  $V$ is motivated by the disjoint union $Y\coprod V$
\end{cor}

This is an immediate consequence of the previous corollary.  It can
also be deduced from the blow up sequence \cite[sect. 9]{manin}, and
this works in any characteristic.

\begin{cor}\label{cor:uniruled}
  A uniruled $n$ dimensional variety is motivated by an $n-1$
  dimensional variety. A unirational $n$ dimensional variety is
  motivated by a variety of dimension less than $n-1$.
\end{cor}

\begin{proof}
  If $X$ is a smooth uniruled $n$-fold, then there is a dominant
  rational map $\PP^1\times Y\dashrightarrow X$ with $\dim Y=n-1$.  By
  resolution of singularities, we can find a sequence of blow ups
  $B_N\to \ldots B_1\to\PP^1\times Y$ along smooth centers and a
  surjective morphism $B_N\to X$.  Then $X$ is motivated by $Y\coprod
  C_1\coprod \ldots C_N$, where $C_i$ are centers of the blow ups.
  This has dimension $n-1$, since the centers have dimension at most
  $n-2$.  A unirational variety is dominated by an iterated blow up of
  $\PP^n$.  So it is motivated by the union of the centers.
\end{proof}

\begin{cor}
  If $Y$ is a smooth projective variety with a $\C^*$-action, then $Y$
  is motivated by the fixed point set $Y^{\C^*}$.
\end{cor}

\begin{proof}
  By the Bialynicki-Birula decomposition \cite{bb}, we can decompose
  $Y$ into a union $Y=\cup Y_i$, where $Y_i$ is a affine space bundle
  over a component of the fixed point set.
\end{proof}

From the discussion following theorem~\ref{thm:gillet}, and the
$G$-equivariant form of resolution of singularities, we get.

\begin{lemma}\label{lemma:XmodG}
  Suppose that the action of a finite group $G$ on a smooth variety
  extends to a compactification, and that the quotient $X/G$ exists in
  $Var_\C$. Then $X/G$ is (weakly) motivated by $\V$ if $X$ is.
\end{lemma}

The quotient $X/G$ always exists as an algebraic space. Thus we could
drop the above requirement by extending the above notions to the
category of algebraic spaces. However, we won't need this.

\section{The coniveau filtration}

Let $FPHS$ be the category of filtered polarizable Hodge structures.
This is additive, but not Abelian.  Given objects $(H,L)$ and $(G,L)$,
we have Tate twists: $$L^p(H(c)) = [L^{p+c}H](c)$$
and tensor
products: $$L^p(H\otimes G) = \sum_{i+j=p} L^i H\otimes L^jG$$

We define the {\em level filtration} $\LL^\dt$ on a pure Hodge
structure $H$ to be $\LL^pH=F^p\cap H_\Q$. This is the largest Hodge
substructure of $H$ satisfying $\LL^p\subseteq F^p$.  If $H$ is pure
of weight $m$, it follows that $\LL^p$ is the maximal substructure
with level at most $|m-2p|$.

\begin{lemma}\label{lemma:polar}
  The operation $V\mapsto \LL^pV$ gives rise to an exact endofunctor
  on the category of polarizable Hodge structures.  The functor
  $V\mapsto (V,\LL^\dt V)$ from $PHS\to FPHS$ is compatible with Tate
  twists and products.
\end{lemma}

\begin{proof}
  The operation $H\mapsto \LL^pH$ is easily seen to be an additive
  functor.  In particular, it preserves direct sums. Since $PHS$ is
  semisimple, this forces exactness. The remaining properties are
  straightforward.
\end{proof}

Let $X$ be a smooth projective variety, the {\em coniveau} filtration
is given by
\begin{eqnarray*}
  N^p H^i(X) &=& \sum_{codimY\ge p} ker[H^i(X)\to H^i(X-Y)]\\
  &=& \sum_{codim Y=q\ge p} im[H^{i-2q}(\tilde Y)(-q)\to H^i(X)]
\end{eqnarray*}
where $Y$ ranges over closed subvarieties; in the second expression
$\tilde Y\to Y$ are chosen desingularizations.  Since the level of
$H^{i-2q}(\tilde Y)(-q)$ is bounded by $i-2p$, we have an inclusion
$$N^p H^i(X) \subseteq \LL^pH^i(X)$$
The generalized Hodge conjecture
asserts that equality holds.  This would imply functoriality of the
coniveau filtration.  Fortunately, this can be checked directly. The
following is proven in \cite{ak}:

\begin{thm}\label{thm:functorialityofN}
  The coniveau filtration $N^\dt$ is preserved by pushforwards,
  pullbacks, and products. More precisely;
  \begin{enumerate}
  \item If $f:X\to Y$ is a map of smooth projective varieties of
    dimensions $n$ and $m$ respectively, then
    $$f_*(N^pH^i(X))\subseteq N^p(H^{i+2(m-n)}(Y)(m-n))$$

  \item If $f$ is as above, then $$f^*(N^pH^i(Y))\subseteq N^pH^i(X)$$
  \item $$N^p(H^i(X))\otimes N^q(H^j(Y))\subseteq
    N^{p+q}H^{i+j}(X\times Y)$$
  \end{enumerate}
\end{thm}

\begin{cor}\label{cor:corr}
  The action of a correspondence preserves the coniveau filtration.
\end{cor}

This allows us to define the coniveau filtration of a motive by
$$N^jH^i(X,p,m) = p_*N^j(H^{i+2m}(X)(m))$$

\section{The conjectures}
We work over $\C$. We recall the basic conjectures as conditions on a
fixed smooth projective variety $X$.

\begin{enumerate}
\item[$D(X)$: ] Homological equivalence coincides with numerical
  equivalence on $X$.
\item[$B(X)$: ] For each $i\le \dim X$, there exists an algebraic
  correspondence inducing an isomorphism $$\nu^i: H^{\dim
    X+i}(X,\Q)\stackrel{\sim}{\to} H^{\dim X-i}(X,\Q)$$
\item[$HC(X)$: ] Any Hodge (i.e. rational $(p,p)$) cycle on $X$ is
  algebraic.
\item[$GHC(X)$:] $N^p H^i(X) = \LL^pH^i(X)$ for all $i,p$.
\item[$AC(X)$: ] All Hodge cycles on $X$ are motivated in the
widest sense (i.e. motivated with respect to $SPVar_\C$).
\end{enumerate}

$D$ and $B$ are among Grothendieck's standard conjectures
\cite{groth,kleiman1,kleiman}.  $B$ is called the Lefschetz standard
conjecture.  $HC$ and $GHC$ are the Hodge and generalized Hodge
conjectures respectively.  $AC$ is due to Andr\'e; it sits in between
the Hodge conjecture and Deligne's conjecture \cite{dm} on the
absoluteness of Hodge cycles.  The Hodge conjecture is well known to
be equivalent to the fullness of the embedding $M_{hom}\to PHS$. A
similar interpretation holds for $AC$ in terms of $M_A\to PHS$.  We
have implications $GHC(X)\Rightarrow HC(X)\Rightarrow D(X)$ and
$D(X\times X)\Leftrightarrow B(X)$ \cite{kleiman1,kleiman}. It is
straightforward to extend some of these conjectures to motives.  Given
$M$ in $M_A$, $GHC(M)$ (respectively $HC(M)$) would assert $N^p H^i(M)
= \LL^pH^i(M)$ for all indices (respectively for $i=2p$).  The
formulations of $HC$ by Jannsen \cite[7.9]{jannsen} and $GHC$ by Lewis
\cite[appendix A]{lewis} for a general variety $X$ are equivalent to
$HC(Gr_0[X])$ and $GHC(Gr_0[X])$ respectively.  We should emphasize
that while technically convenient, these extensions to motives and
general varieties are no stronger than the original conjectures.

The following is a repackaging of results of Andr\'e, Grothendieck,
Jannsen and Kleiman.

\begin{thm}
  The following are equivalent:
  \begin{enumerate}
  \item $M_{hom}(X)$ is semisimple and Abelian.
  \item $M_{hom}(X)\to M_{A}(X)$ is an equivalence.
  \item Numerical equivalence coincides with homological equivalence
    on $X$ and all its powers.
  \item The Lefschetz standard conjecture holds for $X$.
  \end{enumerate}
\end{thm}

\begin{proof}
  The equivalence of (3) and (4) is proven in \cite[prop
  5.1]{kleiman}.  (2) implies (1) by \cite[4.2]{andre}.  (1) implies
  (3) by the first step of the proof of \cite[thm 1]{jannsen2}.
  Finally assume (4). Then conjecture $B$ holds for all powers $X^{N}$
  \cite[prop 4.2]{kleiman}.  A motivated cycle modeled on the category
  generated by $X$, is an expression of the form
  $$\gamma=p_{X^n*}(\alpha\cup *\beta)$$
  where $\alpha,\beta\in
  A(X^{n+m})$, and $*$ is the Lefschetz involution \cite{andre}.
  Since $B(X^{n+m})$ holds, $*\beta$ would be algebraic by \cite[prop.
  1.2]{andre} and \cite{kleiman}.  Therefore $\gamma$ would be
  algebraic.  Thus (2) holds.
\end{proof}

We now come to the main point.

\begin{lemma}\label{lemma:motivatinglemma}
  Suppose that $X$ and $Y$ are smooth projective varieties such that
  $Y$ is motivated by $X$. If $X$ and all its powers satisfies one of
  the conjectures ($D$, $B$, $HC$, $GHC$) stated above, then the
  same conjecture holds for $Y$ and all its powers.
If $Y$ is weakly motivated  by $X$ and $AC$ holds for all powers
of $X$, then it holds for all powers of $Y$.
\end{lemma}

\begin{proof}
  Since $Y^n$ is also motivated by $X$, it suffices to prove the
  conjectures hold for $Y$ alone.

  Suppose that $D(X^n)$ holds for all $n$.  Then the motive $[Y]\in
  M_A(X)$ is a direct summand of some ${\Xi}=\oplus X^{n_i}(j_i)$ with
  complement say $Y'$. Given $\gamma\in H(\Xi)$, let us write
  $\gamma_1$ and $\gamma_2$ for its component with respect to the
  decomposition $H(\Xi)= H(Y) \oplus H(Y')$.  Since $M_A(X)$ is
  equivalent to $ M_{hom}(X)$ by the previous theorem, this
  decomposition of $\Xi$ lies in $M_{hom}(X)$, therefore $\gamma_i$
  are both algebraic if and only if $\gamma$ is.  Suppose that
  $\alpha\in H(Y)$ is an algebraic cycle which is numerically
  equivalent to $0$. We can lift it to a class $\beta\in H(\Xi)$ with
  $\beta_1=\alpha$ and $\beta_2=0$. For any other algebraic cycle
  $\gamma$, we have $\gamma\cdot \beta = \gamma_1\cdot \alpha=0$.
  Therefore $\beta$ is numerically trivial, and consequently
  homologically trivial.

  Since the statements $\forall n\, D(X^n)$ and $\forall n\,B(X^n)$
  are equivalent, case $B$ follows from the previous one.

  Suppose that $ HC(X^n)$ or $ GHC(X^n)$ holds for all $n$. We can
  repeat the previous argument to write $H(\Xi)= H(Y) \oplus H(Y')$.
  Any Hodge cycle $\alpha\in H(Y)$ can be lifted to a Hodge cycle on
  $\beta$ on $\Xi$ with $\beta_1=\alpha$ and $\beta_2=0$. Assuming
  $HC(X^n)$, $\beta$ would have to algebraic, and therefore $\alpha$
  is also algebraic. Assuming $GHC(X^n)$, the equality
  $N^pH(\Xi)=\LL^pH(\Xi)$ forces a similar equality for $Y$.

  Finally, suppose that $ AC(X^n)$ holds for all $n$ and
that $Y$ is weakly motivated by $X$. The argument of
  the previous paragraph with ``algebraic'' replaced by ``motivated''
  and within $M_A(X)^{full}$ shows $AC(Y)$.
\end{proof}

For conjecture $B$, see \cite[thm 5.11]{delb} for a refinement.

\begin{cor}
  The Lefschetz standard conjecture holds for any variety motivated by
  a curve or surface.  In particular, it holds for a uniruled
  threefold, a unirational fourfold. 
\end{cor}

\begin{proof}
  The Lefschetz conjecture for a curve or surface follows from the
  Lefschetz $(1,1)$ theorem. Therefore it holds for a power of a curve
  or a surface by \cite[prop 4.3.1]{kleiman}.  The second statement
  follows from corollary \ref{cor:uniruled}.
\end{proof}

We can recover a result of Lieberman that the Lefschetz conjecture
holds for an Abelian variety, since its cohomology is generated by
$H^1$.  We also note that ``most'' varieties are not motivated by
surfaces by lemma \ref{lemma:schoen}.

\begin{cor}
  If $X$ is weakly motivated by an Abelian variety, then $AC$ holds
for $X$ and all its powers.   
\end{cor}

\begin{proof}
  This follows from \cite[thm 0.6.2]{andre}. 
\end{proof}

We can see that the hypothesis holds for  a unirational threefold by 
corollary \ref{cor:uniruled},  or a  smooth projective variety $X$ whose 
cohomology is generated as an
algebra by $H^1(X)$ by lemma \ref{lemma:HXgend}.
 Additional examples, provided by \cite{andre, andre2},
include K3 surfaces and cubic hypersurfaces of dimension at most 6.

\section{Fourier-Mukai transforms}

We return to the case of a general field $k$. As we saw earlier, in
order to prove that a smooth projective variety $Y$ is motivated by
another such variety $X$, it is necessary to find a suitable
correspondence from a sum of powers of $X$ to $Y$. When $Y$ is a
moduli space of objects on $X$, the correspondence can often be
constructed with the help of a Fourier-Mukai transform or something
close to it.  Fix a sheaf $E$ on $X\times Y$ or more generally an
object in the bounded derived category of coherent sheaves $D(X\times
Y)$.  The {\em Fourier-Mukai transform} with kernel $E$ is the exact
(i.e. triangle preserving) functor $\Phi_E:D(X)\to D(Y)$ given by
$$\Phi_E(\F) = \R p_{Y*}(p_X^*\F\otimes E) $$
where $p_X,p_Y$ denote
the projections.  Given $F\in D(Y\times Z)$, the composition
$\Phi_F\circ \Phi_E$ is again a Fourier-Mukai transform: $$F\circ E =
\R p_{XZ_*}(p_{XY}^*E\otimes p_{YZ}^*F) $$
Furthermore, the functor
$\Phi_E$ has left and right adjoints which can also be realized as
Fourier-Mukai transforms. Specifically, if $E^T\in D(Y\times X)$ is
the ``transpose'' of $E$ and $E^{T*}=\R{\mathcal H}om(E, O_{Y\times
  X})$ its dual, then the right adjoint is $\Phi_{E^{T*}\otimes
  \omega_X[\dim X]}$.  Proofs of these facts can be found in
\cite{mukai,orlov}.

Given an object $E$ in $D(X\times Y)$, we can pass to a
$K_0$-correspondence $\chi( E)=\sum(-1)^i h^i(E)\in K_0(X\times Y)$.
The K\"unneth formula implies $\chi(F\circ E) = \chi(F)\circ \chi(E)$.
Next, we construct a functor, which we call the Mukai functor
$\mu:M_{K_0}\to M_{CH}^u$. It is enough to describe this on
$Cor_{K_0}$.  The putative functor $\mu$ sends $X$ in $M_{K_0}$ to
$[X]$.  Given $e\in K_0(X\times Y)$, define $\mu(e) = ch(e)\cdot
\sqrt{td(X\times Y)}$, where $ch$ is the Chern character
$$ch:K_0(\>)\to CH(\>)\otimes \Q,$$
and $$\sqrt{td(\>)}=1
+\frac{c_1(\>)}{4} + \frac{c_1(\>)^2}{96}+\frac{c_2(\>)}{24}+\ldots $$
is the formal square root of the Todd class of the tangent bundle.

\begin{lemma}
  $\mu$ is a functor.
\end{lemma}

\begin{proof}

  Let $\delta:X\to X\times X$ be the diagonal embedding, and $\Delta =
  im(\delta)$. The classes $O_\Delta\in K_0(X\times X)$ and
  $[\Delta]\in H^*(X\times X)$ represents the identity in their
  respective categories.  From standard properties \cite[ex.
  3.2.4]{fulton}, $$\delta^*td(X\times X) = td(X)^2.$$
  Applying the
  Grothendieck-Riemann-Roch theorem \cite[thm 15.2]{fulton} yields
  \begin{eqnarray*}
    ch(\delta_*O_X)td(X\times X) &=& \delta_*(ch(O_X)td(X))\\
    &=& \delta_*(ch(O_X)\delta^*\sqrt{td(X\times X)}\\
    &=& \delta_*(ch(O_X))\sqrt{td(X\times X)}\\
    &=& [\Delta]\sqrt{td(X\times X)}\\
  \end{eqnarray*}
  Thus $\mu(O_\Delta) = [\Delta]$ as required.

  Given $e\in K_0(X\times Y)$ and $g\in K_0(Y\times Z)$, a second
  application of Grothendieck-Riemann-Roch gives:
  \begin{eqnarray*}
    \mu(g\circ e) &=& ch(p_{XZ*}(p_{XY}^*e\cdot p_{YZ}^*g))\sqrt{td(X\times
      Z)}\\
    &=&p_{XZ*}(ch(p_{XY}^*e\cdot p_{YZ}^*g)\cdot 
    p_X^*\sqrt{td(X)}p_Y^*td(Y)p_Z^*\sqrt{td(Z)})\\
    &=& p_{XZ*}(p_{XY}^*[ch(e)\cdot \sqrt{td(X\times Y)}]
    \cdot p_{YZ}^*[ch(g)\cdot \sqrt{td(Y\times Z)}])\\
    &=& \mu(g)\circ \mu(e)
  \end{eqnarray*}

\end{proof}

The functor $\mu$ is easily seen to be additive. However, it is
not compatible with the tensor structures.
A similar argument involving Grothendieck-Riemann-Roch yields the
following less precise result.

\begin{lemma}\label{lemma:Cherneg}
  Given $e\in K_0(X\times Y)$ and $g\in K_0(Y\times Z)$, the Chern
  classes of $g\circ e$ lie in the algebra generated by
  $\{\epsilon_{i,a}\times \gamma_{j,b}\}$ where $$c_i(e) = \sum
  \epsilon_{ia}\times \epsilon_{ia}'$$
  $$c_i(g) = \sum
  \gamma'_{ia}\times \gamma_{ia}$$
  are the K\"unneth decompositions of
  the above Chern classes.
\end{lemma}

\begin{prop}\label{prop:FourierHC}
  Suppose that $X$ and $Y$ are smooth projective varieties with $E\in
  D(Y\times X)$ an object such that $\Phi_E:D(Y)\to D(X)$ is fully
  faithful.  Then there is a split epimorphism of graded Chow motives
  $$\bigoplus [X](i)^{\oplus n_i} \to [Y].$$
  In particular, $Y$ is
  motivated by $X$.
\end{prop}

For the proof we need.

\begin{lemma}
  Suppose that $F:A\to B$ is a fully faithful functor with a right
  adjoint $G:B\to A$. Then $G\circ F$ is naturally equivalent to the
  identity on $A$.
\end{lemma}

\begin{proof}
  We have $$Hom(M,N) \cong Hom(F(M),F(N)) \cong Hom(M, G\circ F(N))$$
  Thus $N\cong G\circ F(N)$ since they represent the same functor.
\end{proof}

\begin{proof}[Proof of proposition \ref{prop:FourierHC}]
  By the results stated earlier, $\Phi_E$ has a right adjoint of the
  form $\Phi_F$ with $F\in D(Y\times E)$. The previous lemma shows
  that this is a left inverse. Therefore $\mu(f\circ e) = id_Y$, where
  $e=\chi(E)$ and $f=\chi(F)$.  Thus $\mu(f):[X]\to [Y]$ gives a split
  epimorphism in $M_{CH}^u$.  After decomposing $\mu(f)$ into its
  homogeneous components, we get a surjection $\oplus [X](i)^{n_i}\to
  [Y]$ in $M_{CH}$.
\end{proof}

\begin{cor}\label{cor:FourierHC}
  Suppose that $X$ and $Y$ are smooth projective varieties with $E\in
  D(Y\times X)$ an object such that
  \begin{enumerate}
  \item $Ext^i(E_s, E_t)= 0$ for all $i$ when $s\not= t$ (where $E_t=
    E|_{\{t\}\times X}$),
  \item $Hom(E_t, E_t)= k$,
  \item $Ext^i(E_t, E_t)= 0$ for all $i>\dim Y$.
  \end{enumerate}
  Then there is a split epimorphism of graded Chow motives $$\bigoplus
  [X](i)^{\oplus n_i} \to [Y].$$
  In particular, $Y$ is motivated by
  $X$.
\end{cor}

\begin{proof}
  Under the above conditions $\Phi_E$ is fully faithful by a theorem
  of Bondal and Orlov \cite[thm 3.3]{bondal}.
\end{proof}

The hypothesis of the next corollary may seem strange at first glance,
however natural examples of pairs of varieties with equivalent derived
categories exist \cite{bondal, mukai,orlov}.

\begin{cor}
  Suppose that $X$ and $Y$ are smooth projective varieties whose
  derived categories are equivalent as triangulated categories.  Then
  the ungraded Chow motives of $X$ and $Y$ are isomorphic.
  Consequently, $X$ and $Y$ are co-motivated
\end{cor}

\begin{proof}
  We appeal to a theorem of Orlov \cite[thm 3.2.1]{orlov} which shows
  that the equivalence $D(X)\to D(Y)$ and its inverse would be induced
  by Fourier-Mukai transforms.
\end{proof}

The hypothesis of corollary \ref{cor:FourierHC} requires that
$Ext^\dt(E_s,E_t)$ is supported on the diagonal. Unfortunately, this
is rather restrictive. The following alternative form will be applied
later on.

\begin{thm}\label{thm:beauville}
  Let $Y$ and $X$ be smooth projective varieties over a field $k$.
  Let $E\in D(Y\times X)$ be an object such that
  \begin{enumerate}
  \item $Hom(E_s, E_t)= 0$ if $s\not= t$ and $k$ otherwise.
  \item $\dim Ext^1(E_t, E_t) = \dim Y$.
  \item $Ext^i(E_s,E_t) = 0$ for $i>1$.
  \end{enumerate}
  Then $Y$ is motivated by $X$.
\end{thm}

The following proposition occurs implicitly in \cite{beauville}.

\begin{prop}[Beauville]
  Let $Y,X, E$ satisfy above conditions.  Then $[\Delta]= c_{\dim
    Y}(E^*\circ E^T)$ in $CH^*(Y\times Y)$.
\end{prop}

\begin{proof}

  The arguments given in \cite{beauville} carry over with very little
  modification. We set $$F= \R p_{YY*}\R {\mathcal
    H}om(p_{YX}^*E,p_{XY}^*E^T)\cong E\circ E^T $$
  By our assumptions,
  $F$ as above can be represented by a complex of vector bundles
  $f:F^0\to F^1$. For any $(s,t)\in Y\times Y$, we have $$0\to
  Hom(E_s,E_t)\to F^0_{s,t}\stackrel{f(s,t)}{\to} F^1_{s,t}\to
  Ext^1(E_s,E_t)\to 0$$
  The $Hom$ above is supported on the diagonal
  $\Delta$. Thus $\Delta$ can be identified with the degeneracy locus
  of the map $f$.  We note that by our assumptions, the codimension of
  $\Delta$ is $\dim Ext^1(E_t,E_t) = rank F^1- rank F^0+1$. This is
  the expected codimension, therefore we are in a position to compute
  the class $[\Delta]$ by Porteous' formula \cite[thm 14.4]{fulton},
  to obtain formula of the proposition.
\end{proof}

\begin{proof}[Proof of theorem \ref{thm:beauville}]
  This is an immediate consequence of the last proposition, lemma
  \ref{lemma:diagonal} and lemma \ref{lemma:Cherneg}

\end{proof}

\section{GHC for general Jacobians}

We make a short digression to prove the generalized Hodge conjecture
for powers of a general curve. The result may be known to experts, but
we give the proof for lack of a suitable reference. Given a complex
Abelian variety $X$, let $Hdg(X)$ denote the Hodge (or special
Mumford-Tate) group of $H=H^1(X)$. This is the smallest $\Q$-algebraic
subgroup of $GL(H)$ whose real points contain the image of the action
$U(1)\to GL(H\otimes \R)$ induced by the Hodge structure.  Given a
polarization $\psi$ of $X$, the Lefschetz group $Lef(X)$, is the
centralizer of $End(X)\otimes \Q$ in $Sp(H^1(X),\psi)$.  The Lefschetz
group turns out to be independent of the polarization, and it always
contains the Hodge group.  The significance of these groups stems from
the fact that the invariants of $H^*(X^n)$ under $Hdg(X)$
(respectively $Lef(X)$) are precisely the Hodge classes (respectively
sums of products of divisor classes). In particular, $HC(X^n)$ holds
for all $n$ whenever these groups coincide.  Further discussion along
with references can be found in \cite{gordon, murty}.

The characterization of Mumford-Tate groups \cite[p. 43]{dm} together
with \cite[7.5]{deligne} (see also \cite[2.2-2.3]{schoen}) yields:

\begin{lemma}\label{lemma:mono}
  Given a polarized integral variation of Hodge structure $V$ over a
  smooth irreducible complex variety $T$, there exists a countable
  union of proper analytic subvarieties $S\subset T$ such that
  $Hdg(V_t)$ contains a finite index subgroup of the monodromy group
  $$image[\pi_1(S,t)\to GL(V_t)]$$
  for $t\notin S$.
\end{lemma}

\begin{thm}[Hazama] Let $X$ be an abelian variety satisfying $Hdg(X) =
  Lef(X)$ and such that all simple factors are of types I or II in
  Albert's classification, then the generalized Hodge conjecture holds
  for $X$.
\end{thm}

\begin{cor}\label{cor:Hazama}
  If $X$ is as above, then the generalized Hodge conjecture holds for
  all powers of $X$.
\end{cor}

\begin{proof} It can be checked that $Hdg(X^k) = Hdg(X)$.  (This is
  obvious from the Tannakian viewpoint, since $H^1(X)$ and $H^1(X^k) =
  H^1(X)^k$ generate the same tensor category.) Also $Lef(X) =
  Lef(X^k)$ \cite[cor. 4.7]{milne}.  Therefore $X^k$ satisfies the
  same conditions as the theorem.
\end{proof}

\begin{cor}\label{cor:ghc4curves}
  If $E=End(X)\otimes \Q$ is a totally real number field such that
  $dim X/[E:\Q]$ is odd then the generalized Hodge conjecture holds
  for all powers $X$.
\end{cor}

\begin{proof} The conditions imply that $X$ is simple of type I.  The
  equality $Hdg(X) = Lef(X)$ follows from \cite[thm 1]{ribet}.
\end{proof}

\begin{prop}\label{prop:vgjacobian}
  There exists a countable union $S$ of proper Zariski closed sets in
  the moduli space $M_g(\C)$ of curves of genus $g\ge 2$, such that if
  $X\in M_g(\C)-S$ then the generalized Hodge conjecture holds for all
  powers of its Jacobian $J(X)$.
\end{prop}

We shall call such a curve {\em very general}.

\begin{proof}
  Choose $n \ge 3$ and let $M_{g,n}$ be the fine moduli space of
  smooth projective curves of genus $g$ with level $n$ structure
  \cite[13.4]{AO}.  Let $\pi:{\mathcal X}\to M_{g,n}$ be the universal
  curve.  Lemma \ref{lemma:mono} applied to $\R^1\pi_*\Z$ shows that
  there exist a countable union of proper subvarieties $S'\subset
  M_{g,n}(\C)$ such that a finite index subgroup of the monodromy
  group $$\Gamma=image[\pi_1(M_{g,n}, t)\to GL(H^1(\mathcal{X}_t))]$$
  is contained in $Hdg(\mathcal{X}_t)$ for each $t\notin S'$. Let $S$
  be the image of $S'$ in $M_g(\C)$.  By Teichmuller theory, any
  finite index subgroup of $\Gamma$ is seen to be Zariski dense in the
  symplectic group (see \cite[12]{hain}). Hence the Hodge group
  contains the symplectic group whenever $t\notin S$.  But this forces
  $$Hdg(J(\mathcal{X}_t))=Lef(J(\mathcal{X}_t)) =
  Sp(H^1(\mathcal{X}_t)).$$

  Fix $X= \mathcal{X}_t$, with $t$ as above.  We will show that
  $End(J(X))\otimes \Q = \Q$, and this will finish the proof by
  corollary \ref{cor:ghc4curves}.  The natural map $$End(J(X))\otimes
  \Q \to End(H^1(X))$$
  is injective, and the image lies in the ring
  $End_{HS}(H^1(X))$ of endomorphisms of the Hodge structure $H^1(X)$.
  This is contained in the space of $Hdg(X)$-equivariant endomorphisms
  of $H^1(X)$. Since $Hdg(X)$ is the full symplectic group, it acts
  irreducibly on $H^1(X)$.  Therefore Schur's lemma implies that
  $End(X)\otimes \Q = \Q$ as claimed.
\end{proof}

\section{Application to Moduli spaces}

Let $X$ be a smooth projective curve defined over $\C$.  Then the
moduli space of stable vector bundles of coprime rank $n$ and degree
$d$ over $X$ is a smooth projective fine moduli space \cite{seshadri}.
More generally, we can consider the moduli space $M$ of stable
parabolic bundles with respect to a given collection of weights [loc.
cit.].  Under appropriate numerical conditions on $n,d$ and the
weights \cite[sect 2]{BN}, which we assume, $M$ is again a smooth
projective fine moduli space.

\begin{thm}\label{thm:moduliVB}
  With $X$ and $M$ as above, $M$ is motivated by $J(X)$.
\end{thm}

The special case where $M$ is moduli space of vector bundles was due
to del Ba\~no. We give a seperate proof for this case which is
entirely self contained.

\begin{proof}[Proof for vector bundles]
  Since $M$ is fine, there is a Poincar\'e bundle $E$ on $M\times X$.
  This satisfies the hypothesis 
  of theorem~\ref{thm:beauville}, therefore $M$ is motivated by $X$,
  and hence to $J(X)$ by lemma \ref{lemma:XvsAlbX}.
\end{proof}

\begin{proof}[Proof for parabolic bundles]
  Biswas and Raghavendra \cite{BR} have shown that $ H(M)$ is
  generated by the K\"unneth components of Chern classes of certain
  universal sheaves on $X\times M$.  Therefore we can apply lemma
  \ref{lemma:ABtype}.
\end{proof}

The first part of the following is due to Biswas and Narasimhan
\cite{BN}.

\begin{cor}
  $M$ (as above) satisfies the Lefschetz standard conjecture and
 $AC$.
\end{cor}

The following corollaries can be deduced by combining the theorem with
known criteria for the validity of Hodge conjecture for Abelian
varieties.

\begin{cor} If $X$ is
  \begin{enumerate}
  \item a curve of genus 2 or 3,
  \item a curve of prime genus such that the Jacobian is simple, or
  \item a Fermat curve $x^m + y^m +z^m = 0$ with $m$ prime or less
    than $21$, or,
  \item a curve admitting a surjection from a modular curve $X_1(N)$,
  \end{enumerate}
  then the Hodge conjecture holds for $M$. If $X$ is a very general
  curve, then the generalized Hodge conjecture holds for $M$.
\end{cor}

\begin{proof}
  A detailed explanation of the ideas involved can be found in
  \cite{arapura}. In outline, for (3) we can apply a theorem of Shioda
  \cite[thm IV]{shioda}.  The remaining results follow from the
  equality of the Hodge and Lefschetz groups of $J(X)$.  For (2), this
  equality can be obtained from work of Tankeev and Ribet \cite[p
  525]{ribet}. For (1), the equality is due Mumford although
  unpublished.  However, a proof can be found in \cite{moonzar}.  In
  case (4), the equality is given by work of Hazama and Murty
  \cite{hazama1}.

  The last statement follows from proposition~\ref{prop:vgjacobian}.
\end{proof}

Let $X$ be a smooth projective surface. Fogarty has shown that the
Hilbert scheme $M$ of zero dimensional subschemes of fixed length $n$
is smooth and projective (see \cite{gott}).

\begin{thm}[Cataldo-Migliorini]
  $M$ is motivated by $X$.
\end{thm}

\begin{proof}
  Let $X^{(n)} = S^nX$ denote the $n$th symmetric power.  let
  $X^{[n]}=M$ be the Hilbert scheme of zero dimensional subschemes of
  $X$ of length $n$.  There are canonical morphisms $p:X^n\to X^{(n)}$
  and $\psi:X^{[n]}\to X^{(n)}$. The map $\psi$ is birational.  These
  spaces have a natural stratification.  Given a partition $\lambda =
  (n_1,n_2, \ldots n_k)$ of $n$ (i.e. a non strictly decreasing
  sequence of positive integers summing to $n$), let $$X^{(n)}_\lambda
  =\{p(x_1,\ldots x_n)\,|\, x_1= x_2=\ldots = x_{n_1}\not= x_{n_1+1} =
  \ldots = x_{n_1+n_2}\not= \ldots\} $$
  and let $$X^{[n]}_\lambda =
  \psi^{-1} X^{(n)}_\lambda . $$
  These are locally closed subsets of
  $X^{(n)}$ and $X^{[n]}$ which will be regarded as subschemes with
  reduced structure. We will argue that each $X^{[n]}_\lambda$ is
  motivated by $X$.  Then the theorem will follow by corollary
  \ref{cor:strat2}.

  The scheme $X^{[n]}_{(n)}$ parameterizes $0$-dimensional subschemes
  with support at a single point. There is a morphism
  $\pi_n:X^{[n]}_{(n)}\to X$ which sends a subscheme to its support.
  Let $U_k\subset X^k$ be the open subset of $k$-tuples with distinct
  components. For a partition $\lambda= (n_1,\ldots n_k)$ of $n$,
  define $$X_\lambda^{<n>} =U_{k}\times_{X^k}\prod_{i=1}^k
  X^{[n_i]}_{(n_i)}$$
  G\"ottsche \cite[2.1.4, 2.2.4]{gott} has shown
  that $\pi_n$ is a Zariski locally trivial fiber bundle where the
  fiber is smooth, projective and has a cellular decomposition. Then
  corollary~\ref{cor:strat2} implies that $[\prod X^{[n_i]}_{(n_i)}]$
  is motivated by $X$.  $U_k$ is motivated by $X$, since it is the
  complement of a diagonal in $X^k$.  Therefore $X^{<n>}_\lambda$ is
  also motivated by $X$.

  G\"ottsche \cite[2.3.3]{gott} has shown that $ X^{[n]}_\lambda$ is a
  quotient of $X^{<n>}_\lambda$ by a subgroup of $S_n$. It follows
  that $X^{[n]}_\lambda$ is also motivated by $X$ by
  lemma~\ref{lemma:XmodG}.
\end{proof}

\begin{cor}
  The Lefschetz standard conjecture holds for $M$.
\end{cor}

\begin{cor}
  If $X$ is an Abelian surface over $\C$, the Hodge conjecture holds
  for $M$.
\end{cor}

\begin{proof}
  As noted earlier, $HC$ holds for all powers of $X$.
\end{proof}
 
\begin{cor}
  Let $X$ be a smooth projective surface over $\C$ with Kodaira
  dimension $\kappa(X)\le 0$, then the conjecture $AC$ holds for $M$.
\end{cor}

 \begin{proof}
   It suffices by the results of \cite[thm 0.62, 0.63]{andre} to prove
   that $X$, and therefore $M$, is motivated by an Abelian variety, a
   K3 surface or (for trivial reasons) a projective space.  Clearly
   $X$ can be assumed minimal since it is co-motivated with a minimal
   model for it.  Using classification of surfaces
   \cite{beauville-surface}, we see that $X$ rational, ruled over a
   curve $C$, or else there exist a surjective map $S\to X$ with $S$
   Abelian or K3. In the last two cases, $X$ is motivated $J(C)$ or
   $S$ as required.
 \end{proof}

 Let $X$ be an Abelian or K3 surface over $\C$ with an ample line
 bundle $H$. Let $M$ be the moduli space of $H$-stable of rank $r$
 torsion free sheaves with fixed Chern classes $c_1,c_2$. Mukai has
 shown that $M$ is always smooth.  Under appropriate conditions on the
 invariants, $M$ is also projective.  See \cite{huybrechts} for
 further details.

\begin{thm}\label{thm:modsurfaces}
  Let $X$ and $M$ be as in the previous paragraph with $M$ is
  projective.  Then $M$ is motivated by $X$.
\end{thm}

\begin{proof}
  By a theorem of Markman \cite{markman}, $H(M)$ is generated by the
  K\"unneth components of Chern classes of a quasi-universal sheaf
  $\E$ on $X\times M$.  Therefore we can apply lemma
  \ref{lemma:ABtype}.
\end{proof}

\begin{cor}
  If $X$ is Abelian of K3 then $B(M)$ and $AC(M)$ hold, and $HC(M)$
  also holds if $X$ in the Abelian case.
\end{cor}

\end{document}